

%
%
%
%
%
%
%


\magnification=1000
\overfullrule=0pt

\font\eightrm=cmr8

\textfont0=\tenrm
\scriptfont0=\sevenrm
\scriptscriptfont0=\fiverm
\def\rm{\fam0\tenrm}

%
%
\font\teni=cmmi10
\font\seveni=cmmi7
\font\fivei=cmmi5
\textfont1=\teni
\scriptfont1=\seveni
\scriptscriptfont1=\fivei

%

%

\font\eightit=cmti8

%
\font\eightbfti=cmbxti10 at 8pt
\font\bfti=cmbxti10


\font\sc=cmcsc10


%
\font\tenbb=msbm10
\font\sevenbb=msbm7
\font\fivebb=msbm5
\newfam\bbfam
\textfont\bbfam=\tenbb
\scriptfont\bbfam=\sevenbb
\scriptscriptfont\bbfam=\fivebb

%

%
\font\tengoth=eufm10
\font\sevengoth=eufm7
\font\fivegoth=eufm5
\newfam\gothfam
\textfont\gothfam=\tengoth
\scriptfont\gothfam=\sevengoth
\scriptscriptfont\gothfam=\fivegoth

%

%

\font\fourteenbf=cmbx10 at 14pt

\font\eightbf=cmbx8
 

\font\tenbf=cmbx10
\font\sevenbf=cmbx7
\font\fivebf=cmbx5
\newfam\bffam
\textfont\bffam=\tenbf
\scriptfont\bffam=\sevenbf
\scriptscriptfont\bffam=\fivebf
\def\bf{\fam\bffam\tenbf}
%
%



\newif\ifpagetitre           \pagetitretrue
\newtoks\hautpagetitre       \hautpagetitre={\hfil}
\newtoks\baspagetitre        \baspagetitre={\hfil}

\newtoks\auteurcourant       \auteurcourant={\hfil}
\newtoks\titrecourant        \titrecourant={\hfil}

\newtoks\hautpagegauche
         \hautpagegauche={\hfil\the\auteurcourant\hfil}
\newtoks\hautpagedroite
         \hautpagedroite={\hfil\the\titrecourant\hfil}

\newtoks\baspagegauche
         \baspagegauche={\hfil\tenrm\folio\hfil}
\newtoks\baspagedroite
         \baspagedroite={\hfil\tenrm\folio\hfil}

\headline={\ifpagetitre\the\hautpagetitre
\else\ifodd\pageno\the\hautpagedroite
\else\the\hautpagegauche\fi\fi }

\footline={\ifpagetitre\the\baspagetitre
\global\pagetitrefalse
\else\ifodd\pageno\the\baspagedroite
\else\the\baspagegauche\fi\fi }



\def\og{\leavevmode\raise.3ex
     \hbox{$\scriptscriptstyle\langle\!\langle$~}}
\def\fg{\leavevmode\raise.3ex
     \hbox{~$\!\scriptscriptstyle\,\rangle\!\rangle$}}
     


\def\build #1_#2^#3{\mathrel{\mathop{\kern 0pt #1}
     \limits_{#2}^{#3}}}




\def\hfl[#1][#2][#3]#4#5{\kern -#1
 \sumdim=#2 \advance\sumdim by #1 \advance\sumdim by #3
  \smash{\mathop{\hbox to \sumdim {\rightarrowfill}}
   \limits^{\scriptstyle#4}_{\scriptstyle#5}}
    \kern-#3}

\def\vfl[#1][#2][#3]#4#5%
 {\sumdim=#1 \advance\sumdim by #2 \advance\sumdim by #3
  \setbox1=\hbox{$\left\downarrow\vbox to .5\sumdim {}\right.$}
   \setbox1=\hbox{\llap{$\scriptstyle #4$}\box1
    \rlap{$\scriptstyle #5$}}
     \vcenter{\kern -#1\box1\kern -#3}}


\newdimen\sumdim
\def\diagram#1\enddiagram
    {\vcenter{\offinterlineskip
      \def\tvi{\vrule height 10pt depth 10pt width 0pt}
       \halign{&\tvi\kern 5pt\hfil$\displaystyle##$\hfil\kern 5pt
        \crcr #1\crcr}}}







\def\vspace[#1]{\noalign{\vskip #1}}



\def\boxit [#1]#2{\vbox{\hrule\hbox{\vrule
     \vbox spread #1{\vss\hbox spread #1{\hss #2\hss}\vss}%
      \vrule}\hrule}}



\catcode`\@=11
\def\Eqalign #1{\null\,\vcenter{\openup\jot\m@th\ialign{
\strut\hfil$\displaystyle{##}$&$\displaystyle{{}##}$\hfil
&&\quad\strut\hfil$\displaystyle{##}$&$\displaystyle{{}##}$
\hfil\crcr #1\crcr}}\,}
\catcode`\@=12






\vsize= 20.5cm
\hsize=15cm
\hoffset=3mm
\voffset=3mm


\def\tend_#1^#2{\mathrel{
        \mathop{\kern 0pt\hbox to 1cm{\rightarrowfill}}
        \limits_{#1}^{#2}}}




\auteurcourant{\eightbf T.~Hausel, M.~Thaddeus \hfill}


\titrecourant{\hfill\eightbf Mirror partners
arising from integrable systems}


{\parindent=0pt\eightbf


G\'eom\'etrie alg\'ebrique/\hskip -.5mm{\eightbfti Algebraic Geometry}

(Physique math\'ematique/\hskip -.5mm{\eightbfti Mathematical Physics})

\vskip 1.5cm



{\fourteenbf \noindent
Examples of mirror partners arising from integrable systems
}
\parindent=-1mm\footnote{}{
To appear in C. R. Acad.\ Sci.\ Paris.
}

\vskip 15pt


\noindent{\bf Tam\'as HAUSEL~${}^{\hbox{\sevenbf a}}$,
Michael THADDEUS~${}^{\hbox{\sevenbf b}}$}

\vskip 6pt


{\parindent=3mm
\item{${}^{\hbox{\fivebf a}}$}Department of Mathematics, University of
California, Berkeley, CA 94720, USA
\par}

\smallskip

{\parindent=3mm
\item{${}^{\hbox{\fivebf b}}$}Department of Mathematics, Columbia
University, New York, NY 10027, USA
\par}

\medskip


}

\vskip3pt

\line{\hbox to 2cm{}\hrulefill\hbox to 1.5cm{}}

\vskip3pt



{\parindent=2cm\rightskip 1.5cm\baselineskip=9pt \item{\bf
Abstract.~}{\eightrm In this note we present pairs of hyperk\"ahler
orbifolds which satisfy two different versions of mirror symmetry.  On
the one hand, we show that their Hodge numbers (or more precisely,
stringy E-polynomials) are equal.  On the other hand, we show that
they satisfy the prescription of Strominger, Yau, and Zaslow (which in the
present case goes back to Bershadsky, Johansen, Sadov and Vafa): that a
Calabi-Yau and its mirror should fiber over the same real manifold,
with special Lagrangian fibers which are tori dual to each other.  Our
examples arise as moduli spaces of local systems on a curve with
structure group {\eightrm SL(}n\/{\eightrm )}; the mirror is the
corresponding space with structure group {\eightrm PGL(}n\/{\eightrm
)}.  The special Lagrangian tori come from an algebraically
completely integrable Hamiltonian system: the Hitchin
system.~}{\eightrm\copyright~Acad\'emie des Sciences/Elsevier, Paris}
\par}

\vskip 10pt


{\parindent=2cm\rightskip 1.5cm
{\bfti Partenaires miroirs provenants des syst\`emes int\'egrables}
\par}

\vskip 10pt

{\parindent=2cm\rightskip 1.5cm\baselineskip=9pt \item{\bf
R\'esum\'e.~} {\eightit Nous pr\'esentons dans cette note des 
paires de V--vari\'et\'es
hyperk\"ahleriennes satisfaisant deux formulations diff\'erentes
de la sym\'etrie miroir. Nous montrons d'une part que leurs nombres
de Hodge (plus pr\'ecis\'ement, leurs E--polyn\^omes de cordes)
coinc\"{\i}dent. D'autre part, nous montrons qu'elles satisfont
le crit\`ere de Strominger, Yau et Zaslow, c'est--\`a--dire qu'elles
sont fibr\'ees sur la m\^eme vari\'et\'e r\'eelle, de sorte que les
fibres soient des tores Lagrangiens sp\'eciaux duaux entre eux. Nos
exemples se pr\'esentent comme espaces de modules de syst\`emes
locaux sur une courbe avec groupe de structure {\eightrm SL(}n\/{\eightrm)}; 
le miroir est l'espace correspondant avec groupe
de structure {\eightrm PGL(}n\/{\eightrm )}. Les tores Lagrangiens
sp\'eciaux proviennent d'un syst\`eme Hamiltonien compl\`etement
int\'egrable alg\'ebriquement, le syst\`eme de Hitchin.  ~\copyright~Acad\'emie
des Sciences/Elsevier, Paris}
\par}

\vskip1.5pt
\line{\hbox to 2cm{}\hrulefill\hbox to 1.5cm{}}
\medskip

\parindent=3mm
\vskip 5mm


\noindent
{\bf Version Fran\c caise Abr\'eg\'ee}

\bigskip

Il s'agit dans cette note de donner quelques exemples de vari\'et\'es
alg\'ebriques complexes o\`u l'on peut v\'erifier directement les
pr\'edictions de sym\'etrie miroir, au sens de Strominger-Yau-Zaslow
[12].

Soient $M$ et $\hat M$ des $V$-vari\'et\'es Calabi-Yau de dimension
complexe $n$.  On appelle $\hat M$ un {\it partenaire miroir faible} de $M$
s'il existe une $V$-vari\'et\'e $N$ de dimension r\'eelle $n$, et des
applications lisses $\pi: M \to N$, $\hat \pi: \hat M \to N$ telles
que, sur un ouvert dense de $N$, les fibres $L_x = \pi^{-1}(x)$ et
$\hat L_x = \hat \pi^{-1}(x)$ sont des tores Lagrangiens sp\'eciaux,
et que les syst\'emes locaux $\pi_1(L_x)$ et $\pi_1(\hat{L}_x)$ sont en dualit\'e.

Consid\'erons le cas o\`u $M$ est {\it hyperk\"ahler}, c'est-\`a-dire
K\"ahler par rapport \`a trois structures complexes $J_1, J_2, J_3: TM
\to TM$ satisfaisant les relations de commutation
des quaternions imaginaires.

\medskip

{\sc Lemme. -- }{\it
Dans ce cas, $L \subset M$ est Lagrangien sp\'ecial par rapport \`a
$J_1$ s'il est Lagrangien holomorphe par rapport \`a $J_2$. }

\smallskip

Si de plus $M$ est un syst\`eme Hamiltonien compl\`etement
int\'egrable alg\'ebriquement (SHCIA), alors il existe sur un ouvert
dense une fibration par tores Lagrangiens holomorphes, comme on s'y
attend.  L'exemple-cl\'e est le {\it syst\`eme de Hitchin},
que l'on peut d\'ecrire comme suit.

Soit $C$ une courbe lisse complexe projective de genre $g$, $p \in C$
un point de base, et $n$ un entier, le rang.  Pour chaque $d \in {\bf
Z}$, soit $M^d_{\rm Dol}$ l'espace de modules des fibr\'es de Higgs
semistables $(E,\phi)$ sur $C$ satisfaisant $\Lambda^n E \cong {\cal O}(d\,p)$ et
${\rm tr} \, \phi = 0$.  De la meme fa\c con, soit $M^d_{\rm DR}$
l'espace de modules des syst\`emes locaux sur $C$ avec groupe de
structure ${\rm SL(}n{\rm )}$ sur $C \backslash \{ p \}$, dont la
monodromie autour de $p$ est $e^{2 \pi i d/n}$.  Nous renvoyons \`a
Simpson [7,8] pour les d\'efinitions exactes.

\medskip

{\sc Th\'eor\`eme {\rm (Hitchin, Simpson)}. -- }{\it Avec la notation
ci-dessus,

\itemitem{\rm (1)} il existe un hom\'eomorphisme $M^d_{\rm DR} \simeq
  M^d_{\rm Dol}$, et une structure hyperk\"ahlerienne sur la
  partie lisse dont les structures complexes $J_1$ et $J_2$
  sont celles provenant de $M^d_{\rm DR}$ et $M^d_{\rm Dol}$;

\itemitem{\rm (2)} il existe une famille alg\'ebrique equisinguli\`ere
  ${\cal M}^d_{\rm Hod} \to {\bf A}^1$ dont le fibre en z\'ero est
  $M^d_{\rm Dol}$ et dont les autres fibres sont $M^d_{\rm DR}$; il
  existe de plus une
  action de ${\rm GL(1)}$ sur ${\cal M}^d_{\rm Hod}$ qui rel\`eve
  l'action lin\'eaire sur ${\bf A}^1$;

\itemitem{\rm (3)} il existe un morphisme propre $\mu_d: M^d_{\rm Dol}
  \to V$, o\`u $V$ est un espace affine ind\'ependant de $d$, qui fait
  de $M^d_{\rm Dol}$ un SHCIA.

}

\smallskip

\medskip

{\sc Th\'eor\`eme 1. -- }{\it Ce r\'esultat reste valable, sur une
courbe $C$ avec un nombre fini de points choisis, pour les espaces
$M^d_{\rm Dol}$ des fibr\'es de Higgs semistables {\rm paraboliques}
avec des poids fixes, et pour les espaces $M^d_{\rm DR}$ des
syst\`emes locaux {\rm filtr\'es} avec la monodromie correspondante.
}

\smallskip

Nous d\'efinissons les espaces de modules correspondants avec groupe 
de structure ${\rm
PGL(}n{\rm )}$ comme \'etant les quotients les quotients de leurs contreparties
${\rm SL(}n{\rm )}$ par le groupe ab\'elien fini $\Gamma = {\rm
Pic}^0(C)[n] \cong {\bf Z}_n^{2g}$.  Celui-ci op\`ere par
tensorisation, et pr\'eserve toutes les applications qui apparaissent
dans le Th\'eor\`eme plus haut.  Nous indiquerons par un accent
circonflexe le passage au quotient par $\Gamma$: par exemple, le
morphisme $\mu_d: M^d_{\rm Dol} \to V$ d\'escend \`a $\hat \mu_d: \hat
M^d_{\rm Dol} \to V$.

\medskip

{\sc Proposition. -- }{\it Pour $v \in V$ g\'en\'erique, les fibres
$\mu_d^{-1}(v)$ et $\hat \mu_d^{-1}(v)$ sont des espaces principaux
homog\`enes des vari\'et\'es ab\'eliennes polaris\'ees $\mu_0^{-1}(v)$
and $\hat \mu_0^{-1}(v)$, qui sont duales entre elles.  }

\smallskip

Pour un $d$ quelconque, $M^0_{\rm Dol}$ (resp.\ $\hat M^0_{\rm Dol}$)
admet donc une fibration holomorphe par tores duaux \`a ceux de
$\hat M^d_{\rm Dol}$ (resp.\ $M^d_{\rm Dol}$).  En appliquant le
lemme, on trouve que le partenaire miroir de $M^d_{\rm DR}$ (resp.\
$\hat M^d_{\rm DR}$) doit {\^e}tre $\hat M^0_{\rm DR}$ (resp.\ $M^0_{\rm
DR}$).  Puisque ces espaces sont des $V$-vari\'et\'es
hyper\-k\"ahleriennes, \`a cause de la sym\'etrie miroir [8], nous
pouvons pr\'edire que leurs polyn{\^ o}mes de Hodge, ou plutot leurs
$E$-polyn{\^o}mes de cordes [1], not\'es $E_{\rm st}$, doivent etre
\'egaux.  On peut v\'erifier la pr\'ediction dans deux cas, o\`u tous
les espaces qui interviennent sont des $V$-vari\'et\'es.

Consid\'erons d'abord le cas parabolique, en supposant
qu'il existe un point avec structure para\-bo\-lique o\`u les poids
sont suffisament g\'en\'eriques.

\medskip

{\sc Conjecture. -- }{\it
Dans ce cas, pour $c,d \in {\bf Z}$,}
$$
E_{\rm st}(M^c_{\rm DR})
=
E_{\rm st}(M^c_{\rm Dol})
=
E_{\rm st}(\hat M^d_{\rm Dol})
=
E_{\rm st}(\hat M^d_{\rm D})
.$$

\smallskip

\medskip

{\sc Th\'eor\`eme~2. -- }{\it
Ceci est valable pour $n = 2$ et $n = 3$.}

\smallskip

Consid\'erons maintenant l'autre bon cas: celui o\`u il n'y a pas de structure
para\-bo\-lique, mais o\`u le rang $n$ et le d\'egr\'e $d$ sont premiers
entre eux.  Dans ce cas, le partenaire miroir de $M^d_{\rm DR}$,
c'est-\`a-dire $\hat M^0_{\rm DR}$, a de mauvaises singularit\'es, et
la relation miroir n'est pas sym\'etrique.  N\'eanmoins, des
calculs sugg\`erent la conjecture suivante.

Vafa [13] et Vafa-Witten [14] ont montr\'e que, pour un groupe fini
$\Gamma$, un \'el\'ement $\rho \in H^2(\Gamma, {\rm U(1))}$ quelconque
d\'efinit un $E$-polyn{\^o}me de cordes, de fa\c con tordue, d'une
$V$-vari\'et\'e $X/\Gamma$. Ils appellent ce ph\'enom\`ene {\it torsion
discr\`ete}.  Dans notre cas $\Gamma \cong {\bf Z}_n^{2g}$.  Soit $\eta$
le g\'en\'erateur standard de $H^2({\bf Z}_n^2; {\rm U(1)}) \cong {\bf
Z}_n$, et soit $\rho = \sum_{i=1}^g \pi_i^* \eta$ o\`u $\pi_i: \Gamma
\to {\bf Z}_n^2$ la projection sur les facteurs $(2i-1)$'\`eme and
$(2i)$'\`eme.

\medskip

{\sc Conjecture. -- }{\it
Pour $c, d \in {\bf Z}$, }
$$
E_{\rm st}(M^c_{\rm DR})
=
E_{\rm st}(M^c_{\rm Dol})
=
E_{\rm st}^{c \rho}(\hat M^d_{\rm Dol})
=
E_{\rm st}^{c \rho}(\hat M^d_{\rm DR})
.$$

\smallskip

\medskip

{\sc Th\'eor\`eme~3. -- }{\it
Ceci est valable pour $n = 2$ et $n = 3$ si $c$ et $n$, et $d$ et $n$,
sont premiers entre eux.}

\smallskip


\vskip 6pt
\centerline{\hbox to 2cm{\hrulefill}}
\vskip 9pt



The object of this note is to exhibit some striking mathematical 
evidence in support of the Strominger-Yau-Zaslow (SYZ) proposal for mirror
symmetry in string theory [12], which however does not yet have a 
completely satisfactory
precise mathematical formulation.  In this note 
we describe some hyperk\"ahler varieties which
satisfy a mathematical version of some of the conditions laid down by 
these authors to be mirror
partners, and whose stringy mixed Hodge numbers turn out to be equal. The
latter could be considered as a hyperk\"ahler analogue of the
 ``topological mirror test'', which was the first
of the mathematically
surprising findings of physicists pursuing the original mirror symmetry for
Calabi-Yau $3$-folds. The examples we present here also involve a complex 
reductive group $G$ and its Langlands
dual group $\hat G$.  We hope $G$ could be taken to be any reductive group, but
at one point we will need to assume that $G$ is $SL(2)$ or $SL(3)$.

An $n$-dimensional K\"ahler orbifold is said to be {\it Calabi-Yau} if
it is Ricci-flat and it is equipped with a nowhere vanishing
holomorphic $n$-form $\Omega$.  A submanifold $L$ is {\it special
Lagrangian} if it is Lagrangian and ${\rm Im}\, \Omega |_L = 0$.

Let $M$ and $\hat M$ be Calabi-Yau orbifolds of complex dimension $n$.
In this note we call them {\it weak SYZ mirror partners} if there exists an
orbifold $N$ of real dimension $n$ and smooth maps $\pi: M \to N$,
$\hat \pi: \hat M \to N$ so that, if $x$ is a regular value of $\pi$
and $\hat \pi$, $L_x = \pi^{-1}(x)$ and $\hat L_x = \hat \pi^{-1}(x)$
are special Lagrangian tori, and if the local systems of the $\pi_1(L_x)$ and
$\pi_1(\hat{L}_x)$ are in duality, which is equivalent to say that 
there is a smoothly varying diffeomorphism $\hat L^*_x \simeq L^{**}_x$.  
Here for any torus $L$
we define $L^* = {\rm Hom}\,(\pi_1(L),{\rm U(1)})$; the double dual is
used to provide a distinguished basepoint.  If in addition there are
special Lagrangian sections $s$ and $\hat s$ of $\pi$ and $\hat \pi$,
then the basepoints are automatically provided, so that $L_x =
L^{**}_x \simeq \hat L^*_x$ and $\hat L_x = \hat L^{**}_x \simeq
L^*_x$; in this case we call $M$ and $\hat M$ {\it  SYZ mirror
partners}. We note that the original SYZ picture in [12] has extra
conditions on the restricted metrics of these tori, but the metric they 
consider arises in a ``large complex limit'', which does not yet have
a completely satisfactory mathematical formulation. Although we hope that 
our examples will also pass some appropriate metric conditions, for now
we omit this, leaving 
the above fibrewise duality property rather weak.  
Note however that in our hyperk\"ahler case we
find something stronger: the tori are principal homogeneous spaces of abelian 
varieties in duality.  

In general it is a hard analytical problem to find special Lagrangian
tori on $M$.  But in one case an easy lemma shows how to reduce it to
a question of holomorphic geometry.  Suppose that the Riemannian
orbifold $M^{4k}$ is {\it hyperk\"ahler}, that is, it is K\"ahler with the three
K\"ahler forms $\omega_1,\omega_2,\omega_3$ with 
respect to three complex structures $J_1, J_2, J_3: TM \to TM$
satisfying the commutation relations of the imaginary quaternions. Then it
is easy to check that the $2$-form $\omega^c_1=\omega_2+i\omega_3$ is a holomorphic
symplectic form in the complex structure $J_1$. We denote $\Omega_1=(\omega^c_1)^k$
the holomorphic volume form, which makes $(M,J_1)$ into a Calabi-Yau
manifold. We define $\omega^c_i$ and $\Omega_i$ cyclically.

\medskip

{\sc Lemma. -- }{\it
In this situation, $L \subset M$ is special Lagrangian with
respect to $J_1$ and $\Omega_1$ provided that it is holomorphic Lagrangian with
respect to $J_2$ and $\omega^c_2$. }

\smallskip

If, in addition, $M$ is an algebraically completely integrable
Hamiltonian system (ACIHS), then it has a generic holomorphic
Lagrangian torus fibration exactly as desired.  A key example is
furnished by the so-called {\it Hitchin system}, which can be
described as follows.

Fix a smooth complex projective curve $C$ of genus $g$, a basepoint $p
\in C$, and a positive integer $n$, the rank.  For any $d \in {\bf
Z}$, let $M^d_{\rm Dol}$ be the moduli space of semistable Higgs
bundles $(E, \phi)$ on $C$ satisfying $\Lambda^n E \cong {\cal
O}(d\,p)$ and ${\rm tr} \, \phi = 0$.  (If $d=0$, these conditions
mean that $(E, \phi)$ has structure group ${\rm SL(}n{\rm )}$.)  Also
let $M^d_{\rm DR}$ be the space of local systems with structure group
${\rm SL(}n{\rm )}$ on $C \backslash \{ p \}$, having monodromy $e^{2
\pi i d/n}$ around $p$.  We refer to Simpson [7,8] for the definitions
of Higgs bundles, semistability, and so on.

\medskip

{\sc Theorem {\rm (Hitchin [5,6], Simpson [10,11])}. -- }
{\it With the above notation,

\itemitem{\rm (1)} there is a homeomorphism $M^d_{\rm DR} \simeq
  M^d_{\rm Dol}$ and a hyperk\"ahler structure on the smooth locus so
  that the complex structures $J_1$ and $J_2$ correspond to those
  provided by $M^d_{\rm DR}$ and $M^d_{\rm Dol}$;

\itemitem{\rm (2)} there is an equisingular algebraic family ${\cal
  M}^d_{\rm Hod} \to {\bf A}^1$ whose zero fiber is $M^d_{\rm Dol}$
  and whose other fibers are $M^d_{\rm DR}$, and a ${\rm
  GL(1)}$-action on ${\cal M}^d_{\rm Hod}$ lifting the standard action
  on ${\bf A}^1$;

\itemitem{\rm (3)} there is a proper morphism $\mu_d: M^d_{\rm Dol}
  \to V$, where $V$ is an affine space independent of $d$, making
  $M^d_{\rm Dol}$ into an ACIHS.  This morphism has a holomorphic
  Lagrangian section when $d=0$.

}

\smallskip

\medskip

{\sc Theorem 1. -- }{\it A similar result holds, over a curve $C$ with
a finite number of punctures, for the moduli spaces $M^d_{\rm Dol}$ of
semistable {\rm parabolic} Higgs bundles with fixed weights and
$M^d_{\rm DR}$ of {\rm filtered} local systems with the
corresponding monodromy.  }

\smallskip

In fact, similar spaces, and corresponding Theorems, exist for {\it
principal\/} Higgs bundles with any reductive structure group [11]; the
space $M^0_{\rm Dol}$ is the case of ${\rm SL(}n{\rm )}$. The spaces
$M^d_{\rm Dol}$ for $d\neq 0$ are not only used for the construction of the 
${\rm PGL(}n{\rm )}$  moduli spaces below but also will play a decisive role in our main
result. The only other group needed
here is ${\rm PGL(}n{\rm )}$; we define the corresponding moduli spaces
as the quotients of their ${\rm SL(}n{\rm )}$
counterparts by the finite abelian group $\Gamma = {\rm Pic}^0(C)[n]
\cong {\bf Z}_n^{2g}$.  This acts by tensorization, and preserves all
of the maps appearing in the result of Hitchin and Simpson.  The
degree $d$ corresponds to the class of the ${\rm
PGL(}n{\rm )}$-bundle in $H^2(C, {\bf Z}_n) = {\bf Z}_n$.  We denote
the quotient by $\Gamma$ with a hat: for example, the morphism $\mu_d:
M^d_{\rm Dol} \to V$ descends to $\hat \mu_d: \hat M^d_{\rm Dol} \to
V$. 

\medskip

{\sc Proposition. -- }{\it For generic $v \in V$, the fibers
$\mu_d^{-1}(v)$ and $\hat \mu_d^{-1}(v)$ are principal homogeneous
spaces of the polarized abelian varieties $\mu_0^{-1}(v)$ and $\hat
\mu_0^{-1}(v)$.  Indeed, for any $c,d$ there are natural
isomorphisms
$${\rm Pic}^d(\mu_c^{-1}(v)) = \hat \mu_d^{-1}(v)\phantom{.}$$
and
$${\rm Pic}^d(\hat \mu_c^{-1}(v)) = \mu_d^{-1}(v),$$
where ${\rm Pic}^d$ denotes the space of line bundles with first Chern class 
$d$ times the polarization of the underlying polarized abelian variety.}

\smallskip

The relevant abelian varieties can be constructed directly,
using the method of spectral covers of Beauville-Narasimhan-Ramanan
[2] and Hitchin [4], and can be seen to be the generalized Prym variety
of the spectral cover and its dual.

Hence for any $d$, $M^0_{\rm Dol}$ (resp.\ $\hat M^0_{\rm Dol}$) has a
holomorphic torus fibration dual to that of $\hat M^d_{\rm Dol}$
(resp.\ $M^d_{\rm Dol}$).  Using the lemma, one finds that $M^0_{\rm
DR}$ and $\hat M^0_{\rm DR}$, with the corresponding 
holomorphic volume forms,  are SYZ mirror partners, and that
$M^c_{\rm DR}$ and $\hat M^d_{\rm DR}$ are weak SYZ mirror partners
for any $c$ and $d$.

If at any puncture the parabolic weights are distinct, then $\mu_d$
and $\hat \mu_d$ have holomorphic Lagrangian sections, so the
dependence on $c$ and $d$ becomes unimportant.

For a compact hyperk\"ahler manifold $M$, the Hodge numbers are easily
seen to be self-mirror in that $h^{p,q}(M) = h^{n-p,q}(M)$.  One would
therefore hope that the mirror $\hat M$ satisfies $h^{p,q}(M) =
h^{p,q}(\hat M)$.  We will see that (with a suitable notion of Hodge
numbers) the latter equality holds for many of our
mirror partners, even though they are not compact or
self-mirror in the sense above.  We have no explanation as to why this equality holds,
rather than the usual mirror relationship $h^{p,q}(M) = h^{n-p,q}(\hat
M)$.

For spaces which are non-compact and have quotient singularities, the
right notion of Hodge polynomial seems to be the {\it stringy
$E$-polynomial} of Vafa, Zaslow, and Batyrev-Dais [1].  This is a
polynomial in two variables $x$ and $y$, and will be denoted $E_{\rm
st}$.

For the present, we assume either (a) that the weights are
sufficiently generic in the parabolic case, or (b) that $n$ is coprime
to both $c$ and $d$ in the non-parabolic case.  In these cases the
moduli spaces without hats are smooth, so their $E_{\rm
st}$-polynomials reduce to the ordinary $E$-polynomials defined in
terms of Deligne's mixed Hodge structures on the compactly supported
cohomology:
$$E(M) = \sum_{k,p,q} (-1)^{k-p-q} h^{p,q}(H^k_{\rm
cpt}(M))x^p y^q.$$
The moduli spaces with hats are global quotients of the
smooth spaces by $\Gamma$, so their $E_{\rm st}$-polynomials can be defined 
by the orbifold method, according to which
$$E_{\rm st}(X/\Gamma)
= \sum_{\{ \gamma \}} E(X^\gamma/C(\gamma))(xy)^{F(\gamma)}.
\eqno(1)
$$ Here the sum runs over conjugacy classes of $\gamma \in \Gamma$,
$X^\gamma$ is the fixed-point set of $\gamma$ on the smooth variety $X$,
$C(\gamma)$ is the centralizer of $\gamma$, and $F(\gamma)$ is an
integer called the {\it fermionic shift}.  See Batyrev-Dais [1] for
details.

Consider first case (a), when $C$ has at least one puncture where the
weights are sufficiently generic.

\medskip

{\sc Conjecture. -- }{\it
In this case, for any $c,d \in {\bf Z}$,}
$$
E_{\rm st}(M^c_{\rm DR})
=
E_{\rm st}(M^c_{\rm Dol})
=
E_{\rm st}(\hat M^d_{\rm Dol})
=
E_{\rm st}(\hat M^d_{\rm DR})
.$$

\smallskip

\medskip

{\sc Theorem~2. -- }{\it
This is true for $n = 2$ and $n = 3$. }

\smallskip

{\it Sketch of proof.}~--~The proof of the first and last equalities
uses the family ${\cal M}_{\rm Hod} \to {\bf A}^1$ described above,
together with a fiberwise compactification constructed as a quotient
$({\cal M}_{\rm Hod} \times {\bf P}^1)/{\rm GL(1)}$.  This part of the
argument is in fact valid for any $n$ and $d$.

The proof of the middle equality is by explicit calculation.  It is
convenient to subtract $E(\hat M^c_{\rm Dol})$ from both sides.  On
the left, what remains is the contribution of cohomology on $M^c_{\rm
Dol}$ {\it not\/} invariant under the $\Gamma$-action.  This can be
computed from the Bia\l ynicki-Birula stratification, using the
methods of Hitchin [5] to describe the fixed points of the ${\rm
GL(1)}$-action.  On the right, the remaining terms are simply those in
equation (1) above with $\gamma \neq 1$.  These can be calculated
using the methods of Narasimhan-Ramanan [9] to describe the fixed
points of the $\Gamma$-action.  The construction again uses spectral
covers, but since they are unramified it has quite a different flavor
from before.

In the rank 2 case, one finds that
$$E_{\rm st}(M^d_{\rm Dol}) - E(\hat M^d_{\rm Dol})
= E_{\rm st}(\hat M^d_{\rm Dol}) - E(\hat M^d_{\rm Dol})
= 2^{m-1}(2^{2g}- 1)(xy)^{3g-3+m}(1+x)^{g-1}(1+y)^{g-1},$$
where $m$ is the number of punctures whose two weights are distinct.

In the rank 3 case, the formula is similar but more complicated.

\medskip

Now turn to case (b), where there are no punctures, but the rank $n$
is coprime to the degrees $c$ and $d$. In this case $M^c_{\rm DR}$ and
$\hat M^d_{\rm DR}$ are merely weak mirror partners, and the stringy
$E$-polynomials are different.  Hitchin has suggested [7] that weak
mirror partners should be understood as being equipped with
``$B$-fields,'' that is, flat gerbes on $M$ and $\hat M$.  Motivated by this,
we expect that the right stringy Hodge numbers should somehow be twisted with these gerbes. We
suggest that these twistings are occuring due to ``discrete torsion'' as
sketched below. 

Vafa [13] and Vafa-Witten [14] have shown that, for a finite group
$\Gamma$, any element $\rho \in H^2(\Gamma, {\rm U(1))}$ can be used
to define stringy $E$-polynomials, in a twisted sense, of an orbifold
$X/\Gamma$. They refer to this as turning on the {\it discrete
torsion}.  In fact, $\rho$ induces homo\-morphisms $C(\gamma) \to
{\rm U(1)}$ for each $\gamma \in \Gamma$, which give rise to local
systems on $X^\gamma/C(\gamma)$; the $E_{\rm st}$-polynomial of
$X/\Gamma$ then generalizes to a polynomial $E^\rho_{\rm
st}(X/\Gamma)$ where cohomology with coefficients in these local
systems is used in equation (1).

In the present case $\Gamma \cong {\bf Z}_n^{2g}$.
Let $\eta$ be the standard generator of $H^2({\bf Z}_n^2; {\rm U(1)})
\cong {\bf Z}_n$, and let $\rho = \sum_{i=1}^g \pi_i^* \eta$ where
$\pi_i: \Gamma \to {\bf Z}_n^2$ is projection on the $(2i-1)$th and
$(2i)$th factors.  This is in some sense the natural class corresponding to
the symplectic form on $C$.
\medskip

{\sc Conjecture. -- }{\it
For any $c, d \in {\bf Z}$, }
$$
E_{\rm st}(M^c_{\rm DR})
=
E_{\rm st}(M^c_{\rm Dol})
=
E_{\rm st}^{c \rho}(\hat M^d_{\rm Dol})
=
E_{\rm st}^{c \rho}(\hat M^d_{\rm DR})
.$$

\smallskip

\medskip

{\sc Theorem~3. -- }{\it
This is true for $n = 2$ and $n = 3$ when $c$ and $d$ are coprime to
$n$. }

\smallskip

The proof is similar to that of Theorem 2.

\noindent {\bf Acknowledgement.} The authors are grateful for Nigel Hitchin for
suggesting the similarity between [4] and [12] in 1996 and for Pierre Deligne for 
numerous useful comments. 

\bigskip


\centerline{\bf References}

\medskip

{\baselineskip=9pt\eightrm
\item{[1]}Batyrev~V.V., Dais~D.,
Strong McKay correspondence, string-theoretic Hodge numbers and mirror
symmetry,
Topology 35 (1996) 901--929.
\smallskip

\item{[2]}Beauville~A., Narasimhan~M.S., Ramanan~S.,
Spectral curves and the generalised theta divisor,
J. Reine Angew.\ Math.\ 398 (1989) 169--179.
\smallskip

\item{[3]} Bershadsky~M., Johansen~A., Sadov~V., Vafa~C., Topological
reduction of 4{D} {S}{Y}{M} to 2{D}, Nuclear Phys. B\ 448 (1995), 166-186

\item{[4]}Hitchin~N.J.,
Stable bundles and integrable systems,
Duke Math.\ J. 54 (1987) 91--114.
\smallskip

\item{[5]}Hitchin~N.J.,
The self-duality equations on a Riemann surface,
Proc.\ London Math.\ Soc.\ (3) 55 (1987) 59--126.
\smallskip

\item{[6]}Hitchin~N.J.,
Lie groups and Teichm\"uller space,
Topology 31 (1992) 449--473.
\smallskip

\item{[7]}Hitchin~N.J.,
Lectures on special Lagrangian submanifolds,
notes from lectures given at the ICTP School on Differential Geometry,
April 1999, preprint {math.DG/9907034}.
\smallskip

\item{[8]}Morrison~D.,
Geometric aspects of mirror symmetry,
preprint {math.AG/0007090}.
\smallskip

\item{[9]}Narasimhan~M.S., Ramanan~S.,
Generalised Prym varieties as fixed points,
J. Indian Math.\ Soc.\ 39 (1975) 1--19.
\smallskip

\item{[10]}Simpson~C.T.,
Moduli of representations of the fundamental group of a smooth
projective variety, I,
Inst.\ Hautes Etudes Sci.\ Publ.\ Math.\ 79 (1994) 47--129.
\smallskip

\item{[11]}Simpson~C.T.,
The Hodge filtration on nonabelian cohomology,
Algebraic geometry --- Santa Cruz 1995, 217--281,
Proc.\ Sympos.\ Pure Math.\ 62, Part 2,
Amer.\ Math.\ Soc., 1997.
\smallskip

\item{[12]}Strominger~A., Yau~S.-T., Zaslow~E.,
Mirror symmetry is {\eightit T}-duality,
Nuclear Phys.\ B 479 (1996) 243--259.
\smallskip

\item{[13]}Vafa~C.,
Modular invariance and discrete torsion on orbifolds,
Nucl\ Phys.\ B 273
(1986) 592--606.
\smallskip

\item{[14]}Vafa~C., Witten~E.,
On orbifolds with discrete torsion,
J. Geom.\ Phys.\ 15 (1995) 189--214.
\smallskip
}

\end